\documentclass[twoside]{amsart}
\usepackage{amsmath,amsthm,amssymb}

\theoremstyle{plain}
\newtheorem{theorem}{Theorem}[section]
\newtheorem{proposition}[theorem]{Proposition}
\newtheorem{lemma}[theorem]{Lemma}
\newtheorem{corollary}[theorem]{Corollary}

\theoremstyle{remark}
\newtheorem{remark}[theorem]{Remark}
\newtheorem{acknowledgments}{Acknowledgments}

\renewcommand{\r}{{\bf r}}
\newcommand{\A}{{\mathbb A}}
\newcommand{\F}{{\mathbb F}}
\newcommand{\Q}{{\mathbb Q}}
\newcommand{\R}{{\mathbb R}}
\newcommand{\Z}{{\mathbb Z}}

\newcommand{\llceil}{\left\lceil}
\newcommand{\rrceil}{\right\rceil}
\newcommand{\llfloor}{\left\lfloor}
\newcommand{\rrfloor}{\right\rfloor}
\def\bar{\overline}

\newcommand{\cL}{ {\mathcal L} }
\newcommand{\cM}{{\mathcal M}}

\newcommand{\cU}{ {\mathcal U} }
\newcommand{\cV}{ {\mathcal V} }
\newcommand{\cW}{ {\mathcal W} }

\newcommand{\ord}{{\rm ord}}
\newcommand{\sgn}{{\rm sgn} }
\newcommand{\NP }{ {\rm  NP} }     
\newcommand{\HP }{ {\rm  HP} }     
\newcommand{\FP }{ {\rm  FP} }     
\newcommand{\GNP }{ {\rm GNP} }     
\newcommand{\Tr}{ {\rm Tr} }
\newcommand{\Zeta}{{\rm Zeta}}
\newcommand{\N}{{\rm N}}

\begin{document}

\title{$p$-adic variation of $L$ functions of
one variable exponential sums, I
}
\author{Hui June Zhu}
\address{
Hui June Zhu,
Department of mathematics and statistics,
McMaster University,
Hamilton, ON L8S 4K1,
Canada.
}
\email{zhu@cal.berkeley.edu}

\date{October 15 2002}
\subjclass{11,14}

\begin{abstract}
  For a polynomial $f(x)$ in $(\Z_p\cap\Q)[x]$ of degree $d\geq 3$ let
  $L(f\otimes\F_p;T)$ be the $L$ function of the exponential sum of
  $f\bmod p$.  Let $\NP(f\otimes\F_p)$ denote the Newton polygon of
  $L(f\otimes\F_p;T)$.  Let $\HP(\A^d)$ denote the Hodge polygon of
  $\A^d$, which is the lower convex hull in $\R^2$ of the points
  $(n,\frac{n(n+1)}{2d})$ for $0\leq n\leq d-1$.  Let $\A^d$ be the
  space of degree-$d$ monic polynomials parameterized by their
  coefficients.  Let
  $\GNP(\A^d;\F_p):=\inf_{\bar{f}\in\A^d(\F_p)}\NP(\bar{f})$ be the
  lowest Newton polygon over $\F_p$ if exists.  We prove that for
  $p$ large enough $\GNP(\A^d;\F_p)$ exists and we give an explicit
  formula for it.

We prove that there is a Zariski dense open subset $\cU$ defined over
$\Q$ in $\A^d$ such that for $f\in\cU(\Q)$ and for $p$ large enough we
have $\NP(f\otimes\F_p)=\GNP(\A^d;\F_p)$; furthermore, as $p$ goes to
infinity their limit exists and is equal to $\HP(\A^d)$.

Finally we prove analogous results for the space of polynomials
$f(x)=x^d+ax$ with one parameter. In particular, for any nonzero
$a\in\Q$ we show that
$\lim_{p\rightarrow\infty}\NP((x^d+ax)\otimes\F_p)=\HP(\A^d)$.
\end{abstract}

\maketitle

\section{Introduction}\label{S:1}

In this paper $d$ is an integer $\geq 3$.  Let $\A^d$ be the
$d$-dimensional affine space identified with the space of degree-$d$
monic polynomials parameterized by their coefficients.  We always
assume that $p$ is a prime coprime to $d$.
Let $\bar\Q_p$ and $\bar\Z_p$ be the algebraic closure of
$\Q_p$ and its ring of integers respectively.
Let $f(x)$ a polynomial
of one variable in $\A^d(\Z_p\cap\Q)$.  Let
$E(x)=\exp(\sum_{j=0}^{\infty}\frac{x^{p^j}}{p^j})$ be the Artin-Hasse
exponential function. Let $\gamma$ be a root of $\log(E(x))$ in
$\bar\Q_p$ with $\ord_p \gamma=\frac{1}{p-1}$. Then $E(\gamma)$ is a
primitive $p$-th root of unity. Denote it by $\zeta_p$. It is observed
that $\Z_p[\gamma]=\Z_p[\zeta_p]$.  For every $\ell\geq 1$, recall the
exponential sums of the reduction $f\otimes\F_p$ of $f$ modulo $p$
$$S_\ell(f\otimes\F_p):=
\sum_{x\in\F_{p^\ell}}\zeta_p^{\Tr_{\F_{p^\ell}/\F_p}(f(x)\otimes\F_p )}.$$
The $L$-function of the exponential sum
of $f\otimes\F_p$ is defined by
\begin{eqnarray}\label{E:Lfunction}
L(f\otimes\F_p;T):
=\exp(\sum_{\ell=1}^{\infty}S_\ell(f\otimes\F_p)\frac{T^\ell}{\ell}).
\end{eqnarray}
It is well-known (or simply using the Weil Conjecture for curves
combined with (\ref{E:Zeta}) below) that
\begin{eqnarray}\label{E:20}
L(f\otimes\F_p;T)&=&1+b_1T+b_2T^2+\ldots+b_{d-1}T^{d-1}\in\Z[\zeta_p][T].
\end{eqnarray}

Let $\ord_p (\cdot)$ denote the unique extension of the (additive)
$p$-adic valuation in $\Q_p$ to $\bar\Q_p$.  We also denote by
$\ord_p(\cdot)$ the $p$-adic valuation of the {\em content of a power
  series over $\bar\Z_p$} (see \cite[pages 209 and 181]{Lang} for its
standard definition).  Define the {\em Newton polygon} of the
$L$-function of $f\otimes\F_p$, denoted by $\NP(f\otimes\F_p)$, as the
lower convex hull of the points $(n,\ord_p b_n)$ in $\R^2$ for $0\leq
n\leq d-1$, where we set $b_0=1$.  The {\em Hodge polygon} of $f$,
denoted by $\HP(\A^d)$, is the lower convex hull in $\R^2$ of the
points $(n,\frac{n(n+1)}{2d})$ for $0\leq n\leq d-1$.  It is known
that $\HP(\A^d)$ is a lower bound of $\NP(f\otimes\F_p)$ (see
\cite[Propositions 2.2 or 2.3]{Wan:1}) and that if $p\equiv 1\bmod d$
then $\NP(f\otimes\F_p)=\HP(\A^d)$ for every $f\in\A^d(\Q)$ (see
\cite[(3.11)]{Sperber}).

The main result of this paper is Theorems \ref{T:1}, \ref{T:2} and
\ref{T:3}. Theorem \ref{T:1} was a conjecture of Daqing Wan,
proposed in the following form in the number theory seminar at
Berkeley in the fall of 2000 (see also~\cite[Section 2.5]{Wan} for
developments related to this topic).
This theorem follows from Theorem \ref{T:2}.

\begin{theorem}\label{T:1}
  There is a Zariski dense open subset $\cU$ defined over $\Q$ in
  $\A^d$ such that for all $f(x)\in\cU(\Q)$ we have
  $$\lim_{p\rightarrow\infty}\NP(f\otimes\F_p)=\HP(\A^d).$$
\end{theorem}

\begin{remark}
The case $d=3$ follows from \cite[(3.14)]{Sperber},
and the case $d=4$ is discussed in \cite[Corollary 4.7]{Hong}.
The first slope case is proved recently by an elementary method in
\cite{SZ:3} (see also \cite{SZ:2}).  Results concerning Wan's
conjecture ``over $\bar\Q$'' are forthcoming in \cite{Zhu:2}.
\end{remark}

Theorem \ref{T:1} yields an
answer toward questions (in one variable case)
proposed by Katz which asked how the Newton polygon of
the $L$ function of exponential sums
varies with the prime $p$ (see Katz's questions and
Sperber's example on page 151 of \cite[Chapter 5.1]{Katz:2}).  For more
developments in these directions see \cite{Sperber:3},
\cite{Sperber:4} and \cite{Adolphson-Sperber:2} and their
bibliographies.

Let
$X_f: y^p-y=f(x)\otimes\F_p$ be an Artin-Schreier curve over $\F_p$.
The Newton polygon of $X_f$,
denoted by $\NP(X_f\otimes\F_p)$, is
the $p$-adic Newton polygon of the numerator
of the Zeta function $\Zeta(X_f\otimes\F_p;T)$ of $X_f$ over $\F_p$.
It is well-known that (see, for example, \cite[Section VI, (93)]{Bombieri})
\begin{eqnarray}\label{E:Zeta}
\Zeta(X_f\otimes\F_p;T)
=\frac{\N_{\Q(\zeta_p)/\Q}(L(f\otimes\F_p;T))}{(1-T)(1-pT)},
\end{eqnarray}
where the norm $\N_{\Q(\zeta_p)/\Q}$ being interpreted as the product
of the conjugates of $L(f\otimes\F_p;T)$
in $\Q(\zeta_p)$ over $\Q$, the
automorphism acting trivially on the variable $T$.  Thus
$\NP(f\otimes\F_p)$ is precisely equal to $\NP(X_f\otimes\F_p)$ shrunk
by a factor of $\frac{1}{p-1}$ horizontally and vertically, which is
denoted by $\frac{\NP(X_f\otimes\F_p)}{p-1}$. From these remarks the
following corollary is obvious.

\begin{corollary}\label{C:1}
  There exists a Zariski dense open subset $\cU$ defined over $\Q$ in
  $\A^d$ such that for every $f\in \cU(\Q)$ we have
$$\lim_{p\rightarrow\infty}\frac{\NP(X_f\otimes\F_p)}{p-1} = \HP(\A^d).$$
\end{corollary}

\begin{remark}\label{R:13}
a) The behavior of $\NP(x^d\otimes\F_p)$ is well-understood.
For the reader's convenience we describe them briefly below. Let
$\sigma$ be a permutation in the symmetric group
$S_{d-1}$ such that for every $1\leq n\leq d-1$ we let $\sigma(n)$ be
the least positive residue of $pn\bmod d$. Write $\sigma$ as a product
of disjoint cycles (including 1-cycles).  Let $\sigma_i$ be a
$\ell_i$-cycle in $\sigma$.  Let $\lambda_i:=\sum n/(d\ell_i)$ where
the sum ranges over all $n$ in the standard representation of the
$\ell_i$-cycle $\sigma_i$. Arrange $\sigma_i$ in such an order that
$\lambda_1\leq \lambda_2\leq \cdots$.  For every $\sigma_i$ in
$\sigma$ let the pair $(\lambda_i,\ell_i)$ of rational numbers
represent the line segment of (horizontal) length $\ell_i$ and of
slope $\lambda_i$.  The {\em joint of line segments}
$(\lambda_i,\ell_i)$ is the lower convex hull consisting of line
segment $(\lambda_i,\ell_i)$'s connected at their end-points.  The
eigenvalues of $L(x^d\otimes\F_p;T)$ are Gauss sums (see \cite[chapter
III]{Koblitz:2}). By Stickelberger theorem (see \cite[Chapter
6]{Washington}), the $p$-adic Newton polygon of $L(x^d\otimes\F_p;T)$
and hence $\NP(x^d\otimes\F_p)$ is the joint of $(\lambda_i,\ell_i)$'s.
(I thank Kiran Kedlaya for discussions here.)

b) For every $d\geq 3$ the
Newton polygon $\NP(x^d\otimes\F_p)$ does not have a limit as $p$
approaches $\infty$. Indeed, it is clear from the above that
for $p\equiv 1\bmod d$ the Newton polygon is equal to the Hodge
polygon while for $p\equiv -1\bmod d$ the Newton polygon is a straight
line of slope $1/2$.
\end{remark}

This paper is organized as follows.  In section \ref{S:2} notations
and terminologies are introduced.  Using Dwork's $p$-adic analysis, we
define {\em Fredholm polygon} of $f(x)$ over $\F_p$ and show that it
is equal to $\NP(f\otimes\F_p)$.  Section \ref{S:3} is a key step in
the proof, it constructs an $n$-th {\em generic polynomial}, denoted by
$f_n^{t_n}$, proves that they are nonzero and hence defines some Zariski
dense open subset $\cV_n$ in $\A^{d-1}$. It is recommended that the
reader skips section \ref{S:3} at first reading and continues with
section \ref{S:4}, where we immediately apply Dwork's $p$-adic theory
to determine the Fredholm polygon.  In section \ref{S:Proof} we prove
that in some Zariski dense open subset the Fredholm polygon and Newton
polygon coincide if $p$ is large enough.  We prove Theorem \ref{T:2}
there. Finally in section \ref{S:6} we study the families
$f(x)=x^d+ax$ and prove Theorem \ref{T:3} there.

\begin{acknowledgments}
  It is my great pleasure to thank Alan Adolphson and Steven Sperber
  who exposed me to Dwork theory during the Dwork trimester in Italy
  (2001).  I thank Hanfeng Li, Daqing Wan and the referees for
  careful reading and very helpful comments to earlier versions.  Most of all
  I thank Bjorn Poonen for generously sharing ideas and answering
  questions. This research was partially supported by UC Berkeley and
  a grant of Bjorn Poonen from the David and Lucile Packard Foundation.
\end{acknowledgments}

\section{Dwork p-adic theory}\label{S:2}

The fundamental material in our exposition follows
\cite[Sections II and III]{Bombieri}
(see also \cite{Dwork} \cite{Dwork2}
and \cite{Adolphson-Sperber}).
Recall that $p$ is a prime number coprime to $d$.
Let $\bar f(x)=x^d+\sum_{i=1}^{d-1}\bar a_ix^i
\in\F_p[x]$.
Let $f(x)=x^d+\sum_{i=1}^{d-1}a_ix^i\in (\Z_p\cap\Q)[x]$
and $a_d=1$ such that reduction of $f(x)$
at $p$ is equal to $\bar f(x)$.
For any $a_0\in\Z_p\bigcap\Q$, by a simple computation with
(\ref{E:Lfunction}),
one easily concludes that
$L((f+a_0)\otimes\F_p;T)=L(f\otimes\F_p;\zeta_p^{a_0}T)$.
Thus we have
\begin{eqnarray}\label{E:98}
\NP((f+a_0)\otimes\F_p)=\NP(f\otimes\F_p).
\end{eqnarray}

Write $\vec{\hat{a}}=(\hat{a}_1,\ldots,\hat{a}_{d-1})$ where
$\hat{a}_i$ is the Teichm\"uller lifting of $\bar a_i$, that is,
$\hat{a}_i\equiv a_i\bmod p$ and $\hat{a}_i^p=\hat{a}_i$.
Let $\vec{a}:=(a_1,\ldots,a_{d-1})\in (\Z_p\cap\Q)^{d-1}$.
Let $\theta(x)=E(\gamma x)$, where $E(\cdot)$ and $\gamma$ as
defined in section 1.
Then we may write
$\theta(x)=\sum_{m=0}^{\infty}\lambda_mx^m$
for $\lambda_m\in\Z_p[\zeta_p]$.
Note the following properties,
\begin{eqnarray}\label{E:22}
\ord_p \lambda_m&\geq&\frac{m}{p-1};
\end{eqnarray}
for $0\leq m\leq p-1$ we have,
\begin{eqnarray}\label{E:33}
\lambda_m=\frac{\gamma^m}{m!}\mbox{ and  } \ord_p \lambda_m=\frac{m}{p-1}.
\end{eqnarray}
Let $\vec{A}=(A_1,\ldots,A_{d-1})$ be a vector of variables and
$\vec{m}=(m_1,\ldots,m_{d-1})$.
Write $\vec{A}^{\vec{m}}$ for the monomial
$A_1^{m_1}\cdots A_{d-1}^{m_{d-1}}$.
Let $G_n(\vec{A})=0$ for $n<0$.
For every integer $n\geq 0$ let
\begin{eqnarray}\label{E:3}
G_n(\vec{A}):=\sum_{\substack{m_\ell\geq 0\\\sum_{\ell=1}^{d}\ell m_\ell=n}}
\lambda_{m_1}\cdots \lambda_{m_d}\vec{A}^{\vec{m}}.
\end{eqnarray}
Clearly we observe
that $G_n(\vec{A})\in \Z_p[\zeta_p][\vec{A}]$,
that is, $G_n(\vec{A})$ is a polynomial in variable $\vec{A}$
and with coefficients in $\Z_p[\zeta_p]$.
For all integers
$m_1,\ldots,m_d\geq 0$ such that $\sum_{\ell=1}^{d}\ell m_\ell=n$,
we have $d(m_1+\cdots +m_d)\geq \sum_{\ell=1}^{d}\ell m_\ell = n$
and so
$\min(m_1+\cdots +m_d)=\llceil\frac{n}{d}\rrceil$.
Therefore by (\ref{E:3}) we have
\begin{eqnarray}\label{E:Bound}
\ord_p G_n(\vec{A})
\geq
\frac{\min(m_1+\cdots +m_{d})}{p-1}
\geq
\frac{\llceil\frac{n}{d}\rrceil}{p-1}
\geq \frac{n}{d(p-1)}.
\end{eqnarray}
Let $G(X):=\prod_{i=1}^{d}\theta(\hat a_i X^i)\in\Z_p[\zeta_p][[X]]$.
We have
\begin{eqnarray*}
G(X)=
(\sum_{m_1=0}^{\infty}\lambda_{m_1}\hat a_1^{m_1}X^{m_1})
\cdots
(\sum_{m_d=0}^{\infty}\lambda_{m_d}\hat a_d^{m_d}X^{dm_d})
=\sum_{n=0}^{\infty}G_n(\vec{\hat a})X^n.
\end{eqnarray*}
Let $C_0(\vec{A})=1$, and for every $n\geq 1$ let
\begin{eqnarray}\label{E:1}
C_n(\vec{A}):=\sum_{1\leq u_1<u_2<\ldots<u_n}
\sum_{\sigma\in S_n}
\sgn(\sigma)\prod_{i=1}^{n}G_{pu_i-u_{\sigma(i)}}(\vec{A}),
\end{eqnarray}
where $\sgn(\sigma)$ is the signature of the permutation $\sigma$
in the $n$-th symmetric group $S_n$.
It can be verified that
this definition makes sense and that $C_n(\vec{A})\in
\Z_p[\zeta_p][[\vec{A}]]$.

\begin{lemma}\label{L:Dwork}
Let $p$ be a prime coprime to $d$.
For $f(x)=x^d+\sum_{i=1}^{d-1}a_ix^i\in(\Z_p\cap\Q)[x]$,
write $\vec{a}=(a_1,\ldots,a_{d-1})$. Let
$\vec{\hat{a}}=(\hat a_1,\ldots,\hat a_{d-1})$
be Teichm\"uller lifting of $\vec{\bar a}=(\bar a_1,\ldots,\bar a_{d-1})$.
Then
\begin{eqnarray}
\label{E:Dwork}
L(f\otimes\F_p;T)&=&1+b_1(\vec{a})T+\cdots +b_{d-1}(\vec{a})T^{d-1}
\\
&=&\frac{1+\sum_{n=1}^{\infty}(-1)^nC_n(\vec{\hat a})T^n}
{(1-pT)(1+\sum_{n=1}^{\infty} (-1)^nC_n(\vec{\hat a})p^nT^n)},
\nonumber
\end{eqnarray}
where $b_1(\vec{a}),\ldots,b_{d-1}(\vec{a})\in\Z[\zeta_p]$.
\end{lemma}
\begin{proof}
The first equality is a rephrase of (\ref{E:20}).
For every positive integer $\ell$ let
$$S^*_\ell(f\otimes\F_p):=
\sum_{x\in\F_{p^\ell}^*}\zeta_p^{\Tr_{\F_{p^\ell}/\F_p}(f(x)\otimes\F_p)}.$$
Let
$$
L^*(f\otimes\F_p;T)
:=\exp(\sum_{\ell=1}^{\infty}S^*_\ell(f\otimes\F_p)\frac{T^\ell}{\ell}).
$$
Note that $S^*_\ell(f\otimes\F_p)=S_\ell(f\otimes\F_p)-1$ so
\begin{eqnarray}\label{E:Lstar}
L^*(f\otimes\F_p;T)
&=&\exp(\sum_{\ell=1}^\infty (S_\ell(f\otimes\F_p)-1))\nonumber\\
&=&(1-T)\exp(\sum_{\ell=1}^\infty S_\ell(f\otimes\F_p)\frac{T^\ell}{\ell})
\nonumber\\
&=&(1-T)L(f\otimes\F_p;T).
\end{eqnarray}
For any $c>0$ and $b\in \R$ let
$\cL(c,b)$ be the set of power series defined by
$$\cL(c,b):=\{
\sum_{n=0}^{\infty}A_nX^n\mid A_n\in \Q_p(\zeta_p),
\ord_p  A_n\geq \frac{cn}{d}+b\}.$$
Let $\cL(c):=\bigcup_{b\in\R}\cL(c,b)$.
{}From (\ref{E:Bound}) we have
$G(X)=\sum_{n=0}^{\infty}G_n(\vec{\hat a})X^n$ lie in $\cL(1/(p-1))$.
For any $\sum B_n X^n$ in $\cL(c)$,
let $\psi$ be the Hecke operator from $\cL(c)$ to $\cL(cp)$
given by
$\psi(\sum B_n X^n) = \sum B_{pn}X^n$.
Let $\alpha_1:=\psi\cdot G(X)$ be the endomorphism of
$\cL(p/(p-1))$ defined by the composition of the multiplication map by
$G(X)$ then $\psi$, namely,
$$
\alpha_1\left(\sum_{i=0}^{\infty}B_iX^i\right)=\sum_{i=0}^{\infty}
\left(\sum_{j=0}^{\infty}G_{pi-j}(\vec{\hat a})B_j \right)X^i.
$$
Choose the standard monomial basis $\{1,x,x^2,\ldots\}$
for the $p$-adic space $\cL(p/(p-1))$.
Then the $\Q_p(\zeta_p)$-endomorphism
$\alpha_1$ of $\cL(p/(p-1))$ has a matrix
representation by $\{G_{pi-j}(\vec{\hat a})\}_{i,j\geq 0}$.
We denote this matrix by $F_1$.
By the Dwork trace formula (see \cite[Section III]{Bombieri})
we have
\begin{eqnarray*}
L^*(f\otimes\F_p;T)&=&\frac{\det(1-F_1T)}{\det(1-F_1pT)}.
\end{eqnarray*}
For the first row (i.e., $i=0$) of $F_1$, we have
$G_{pi-j}(\vec{\hat a})=0$ for all $j\geq 1$
and $G_{0}(\vec{\hat a})=1$.
By (\ref{E:1}) we have
$$\det(1-F_1T)
=(1-T)\det(1-\{G_{pi-j}(\vec{\hat a})T\}_{i,j\geq 1})
=(1-T)\sum_{n=0}^{\infty}(-1)^nC_n(\vec{\hat a})T^n.$$
Therefore, by (\ref{E:Lstar}) we have
\begin{eqnarray*}(1-T)L(f\otimes\F_p;T)
=L^*(f\otimes\F_p;T)=\frac{(1-T)\sum_{n=0}^{\infty}(-1)^nC_n(\vec{\hat
a})T^n}{ (1-pT)\sum_{n=0}^{\infty}(-1)^nC_n(\vec{\hat a}) p^nT^n}.
\end{eqnarray*}
By simplification of the above formula, our assertion follows.
\end{proof}

\begin{proposition}\label{P:Dwork}
Let the {\em Fredholm polygon} of $f\otimes\F_p$,
denoted by $\FP(f\otimes\F_p)$,
be the lower convex hull of points $(n,\ord_pC_n(\vec{\hat{a}}))$
in $\R^2$ for $0\leq n\leq d-1$.
Then $$\NP(f\otimes\F_p)=\FP(f\otimes\F_p).$$
\end{proposition}
\begin{proof}
By (\ref{E:Dwork}) we have
\begin{eqnarray*}\label{E:Dwork2}
L(f\otimes\F_p;T)(1-pT)(1-C_1pT+C_2p^2T^2-\cdots)
& = &1-C_1T+C_2T^2-\cdots.
\end{eqnarray*}
The ($p$-adic) Newton polygon of $1-C_1T+C_2T^2-\cdots$
has only positive slopes (see \cite[III]{Bombieri}),
so the Newton polygon of
$1-C_1pT+C_2p^2T^2-\cdots$ has every slope $>1$.
On the other hand, the Newton polygon of $L(f\otimes\F_p;T)$
is symmetric in the sense that for every slope segment $\alpha$
there is a slope segment $1-\alpha$ of the same
horizontal length. This property is derived from
the same
fact for Newton polygons of Zeta functions of
abelian varieties and hence of Artin-Schreier curves
(see, for example, \cite[Introduction]{LiOort}).
Thus the slopes of $\NP(f\otimes\F_p)$
are positive and $<1$.
Note that the power series $1-C_1T+C_2T^2-\cdots$
is entire (see \cite[page 121]{Koblitz} for a proof),
so are the three factors on the left-hand-side.
By the $p$-adic Weierstrass preparation theorem
(see \cite[IV.4 Theorem 14]{Koblitz}),
$\NP(f\otimes\F_p)$ coincides with the
p-adic Newton polygon of $1-C_1T+\cdots+(-1)^{d-1}C_{d-1}T^{d-1}$.
\end{proof}

We remark that it is {\em not} generally true that
$L(f\otimes\F_p;T)
=1-C_1(\vec{\hat{a}})T+\cdots+(-1)^{d-1}C_{d-1}(\vec{\hat{a}})T^{d-1}$.

\section{Generic polynomials and
Zariski dense subsets}\label{S:3}

The following notations and conventions are adopted for the rest of
the section.  Given a polynomial as a sum (or several sums) of
polynomials, its {\em formal expansion} means
the formal summation
of its monomials (so one does not do {\em ``arithmetic''},
e.g, cancellations, among its terms).
For any
$\vec{m}=(m_1,\ldots,m_{d-1})\in\Z_{\geq 0}^{d-1}$ let
$|\vec{m}|=\sum_{\ell=1}^{d-1}m_\ell$ and $\vec{m}!=m_1!\cdots
m_{d-1}!$. Fix an integer $r$ with $1\leq
r\leq d-1$ and $\gcd(d,r)=1$.  Let $1\leq n\leq d-1$.

\subsection{The residue matrix $\r_n$}\label{S:31}

Let $1\leq i,j\leq d-1$. Let $r_{ij}$ be the least nonnegative
residue of $-(ri-j)\bmod d$. That is,
$r_{ij}:=d\llceil\frac{ri-j}{d}\rrceil-(ri-j)$.
Let $r'_{ij}$ be the least nonnegative residue of $ri-j\bmod d$.

\begin{lemma}\label{L:r_n}
Let $\r_n$ be the matrix $\r_n:=\{r_{ij}\}_{1\leq i,j\leq n}$.
Then $0\leq r_{ij}\leq
  d-1$ and there are no two identical entries in any row (column).  In
  $\r_{d-1}$ for every $1\leq i\leq d-1$ one has $r_{ij}=0$ if and
  only if $j=r'_{i1}+1$.
\end{lemma}
\begin{proof}
By definition, $r_{ij}$ is
the least non-negative residue of $-(ri-j)\bmod d$ so we have $0\leq
r_{ij}\leq d-1$.
We prove for rows.
The argument for columns is almost identical.
Suppose we have $r_{ij}=r_{ij'}$ then
$ri-j\equiv ri-j'\bmod d$ by definition. Then $j\equiv j'\bmod d$.
Since $1\leq j,j'\leq n\leq d-1$ we have $j=j'$. So there are no
identical entries in any row of $\r_n$.
Note that $r$ is coprime to $d$ so for every
$1\leq i\leq d-1$ there is a unique $1\leq j\leq d-1$
(more precisely $j=r'_{i1}+1$)
such that $ri\equiv j\bmod d$.
This is equivalent to
$r_{ij}=0$ by definition. This proves the last assertion.
Moreover,
\end{proof}

Let $A_d$ be an auxiliary variable. Define a
homogeneous auxiliary polynomial
$\dot{D}_n:=\det(\{A_{d-r_{ij}}\}_{1\leq i,j\leq n})$
of degree $n$ in $\Q[A_1,\ldots,A_d]$.
Note that
\begin{equation}
\dot{D}_n=\sum_{\sigma\in S_n}\sgn(\sigma)
                    \prod_{i=1}^{n}A_{d-r_{i,\sigma(i)}}
=\sum_{\sigma\in S_n}\sgn(\sigma)
\prod_{k=0}^{d-1}A_{d-k}^{\#\{1\leq i\leq n|r_{i,\sigma(i)}=k\}}.
\label{E:D_n}
\end{equation}

\begin{lemma}\label{L:Lexi}
There is a unique highest-lexicographic-order-monomial in the
formal expansion of $\dot{D}_n$ in $\Q[A_1,\ldots,A_d]$.
\end{lemma}

\begin{proof}
It is a combinatorial problem and we shall give an intuitive proof.
We shall do so by verifying the correctness of the following
algorithm which can really be used to obtain the desired
highest-lexicographic-order-monomial.

Fix a residue matrix $\r_n$. Let $\sigma$ be a
permutation in $S_n$ awaiting to be defined.
For every entry in $\r_n$ with
$r_{i_0,j_0}=0$, assign $\sigma(i_0):=j_0$
and cross off the $i_0$-row and the $j_0$-column;
Let $\ell_0$ be the number of all such entries.
For every leftover entry in $\r_n$ with
$r_{i_1,j_1}=1$, assign $\sigma(i_1):=j_1$
and cross off the $i_1$-row and the $j_1$-column;
Let $\ell_1$ be the number of all such entries.
Continue this process until all entries are crossed off.

It is straightforward to verify that
this algorithm uniquely defines a permutation
$\sigma$ by the first statement in Lemma \ref{L:r_n}.
Moreover,
$\sigma$ yields the highest-lexicographic-order-monomial.
Indeed, from (\ref{E:D_n})
one notes that $\ell_0$ is the highest-$A_d$-exponent
in the formal expansion of $\dot{D}_n$;
and $\ell_1$ is the highest-$A_{d-1}$-exponent in
a monomial containing $A_d^{\ell_0}$; and so on.
Thus $\sigma$ yields
the (unique) highest-lexicographic-order-monomial
$A_d^{\ell_0}A_{d-1}^{\ell_1}\cdots A_1^{\ell_{d-1}}$
in the formal expansion of $\dot{D}_n$.
\end{proof}


\begin{lemma}\label{L:evaluation}
Let $M$ be the (unique) highest-lexicographic-order-monomial
of formal expansion of $\dot{D}_n$
derived in Lemma \ref{L:Lexi}.
Then $M|_{A_d=1}$ is the (unique)
highest-lexicographic-order-monomial of lowest degree
in the formal expansion of $\dot{D}_n|_{A_d=1}$.
\end{lemma}
\begin{proof}
It is clear that
the evaluation map $\dot{D}_n\rightarrow \dot{D}_n|_{A_d=1}$
(on the formal expansions) yields a bijective map sending the set of
highest-$A_d$-exponents monomials in the formal expansion of
$\dot{D}_n$ to the set of
lowest-degree-monomials in the formal
expansion of $\dot{D}_n|_{A_d=1}$. Applying the same argument for
the rest of the variables inductively, we conclude our assertion
immediately.
\end{proof}

\subsection{The $n$-th generic polynomial $f^{t_n}_n$}

For any $0\leq s\leq n$ one obtains a
nonempty subset in $\Z_{\geq 0}^{d-1}$
\begin{eqnarray*}
\cM^s_{ij}&:=&\{\vec{m}=(m_1,m_2,\ldots,m_{d-1})\in\Z_{\geq 0}^{d-1}\mid
\sum_{\ell=1}^{d-1}\ell m_{d-\ell}=r_{ij}+ds\}.
\end{eqnarray*}
Recall $\vec{A}:=(A_1,\ldots,A_{d-1})$, and
$\vec{A}^{\vec{m}}:=A_1^{m_1}\cdots A_{d-1}^{m_{d-1}}$.
For $1\leq i,j\leq d-1$ let
\begin{equation}\label{E:delta}
\delta_{ij}:=\left\{
\begin{array}{ll}
0&\mbox{for $j< r'_{i1}+1$}\\
1&\mbox{for $j\geq r'_{i1}+1$}.
\end{array}
\right.
\end{equation}
For $0\leq s\leq n$ and $1\leq i,j\leq n$ define an auxiliary polynomial
\begin{eqnarray}
\label{E:H}
H^s_{ij}(\vec{A})&:=& \sum_{\vec{m}\in\cM^s_{ij}}
h^s_{\vec{m},i,j}
\vec{A}^{\vec{m}}
\end{eqnarray}
where
$$h^s_{\vec{m},i,j}:=
\frac{(\frac{r_{i1}-1}{d}+n)(\frac{r_{i1}-1}{d}+n-1)\cdots
  (\frac{r_{i1}-1}{d}-\delta_{ij}+s+1-|\vec{m}|)}{\vec{m}!}.$$

\begin{lemma}\label{L:H}
Let $0\leq s\leq n$ and $1\leq i,j\leq n$.

a). The polynomial $H^s_{ij}(\vec{A})$ in $\Q[\vec{A}]$
is nonzero and is supported on every $\vec{m}\in
\cM^{s}_{ij}$. The degree of its monomials ranges from
$s+\llceil\frac{r_{ij}+s}{d-1}\rrceil$ up to
$r_{ij}+ds$, where the maximal degree is attained at exactly one
monomial $A_{d-1}^{r_{ij}+ds}$ while the minimal degree is attained
at one or more monomials.

b). The polynomial $H^s_{ij}(\vec{A})$ has a
constant term if and only if $s=r_{ij}=0$;
it has a linear term
if and only if $s=0$ and $r_{ij}\neq 0$, in which
case this linear monomial is exactly $A_{d-r_{ij}}$. \\
\end{lemma}
\begin{proof}
a). Since $\gcd(r,d)=1$ we have $-(ir-1)\not\equiv 1\bmod d$.  Hence
$\frac{r_{i1}-1}{d}\not\in\Z$ and $h^s_{\vec{m},i,j}\neq 0$.
Now it remains to show
\begin{eqnarray*}
\max_{\vec{m}\in \cM^s_{ij}}|\vec{m}|=r_{ij}+ds, &&
\min_{\vec{m}\in \cM^s_{ij}}|\vec{m}|=s+\llceil\frac{r_{ij}+s}{d-1}\rrceil\\
\end{eqnarray*}
For $\vec{m}\in \cM^{s}_{ij}$ we have
$|\vec{m}|\leq \sum_{\ell=1}^{d-1}\ell m_{d-\ell}=r_{ij}+ds$
and the equality holds precisely for $m_1=\cdots=m_{d-2}=0$ and
$m_{d-1}=r_{ij}+ds$.
For $\vec{m}\in\cM^{s}_{ij}$ one has clearly $(d-1)|\vec{m}|\geq
\sum_{\ell=1}^{d-1}\ell m_{d-\ell}=r_{ij}+ds$.
So $$|\vec{m}|\geq
\llceil\frac{r_{ij}+ds}{d-1}\rrceil =s+\llceil\frac{r_{ij}+s}{d-1}
\rrceil.$$
It is easy to see that there are
$m_1,\ldots, m_{d-1}\geq 0$ satisfying
\begin{eqnarray}
\label{E:combine}
\sum_{\ell=1}^{d-1} m_{\ell}=\llceil\frac{r_{ij}+ds}{d-1}\rrceil
\quad \mbox{and}\quad
\sum_{\ell=1}^{d-1}(d-\ell)m_\ell=r_{ij}+ds.
\end{eqnarray}
For example, let $\kappa$ be the least non-negative
residue of $-(r_{ij}+ds)\bmod (d-1)$ then
let $m_1=\llceil\frac{r_{ij}+ds}{d-1}\rrceil-\kappa$,
$m_2=\kappa$ and let the rest $m_\ell=0$.
This says that there are $\vec{m}\in \cM^{s}_{ij}$
with $|\vec{m}|=\llceil\frac{r_{ij}+ds}{d-1}\rrceil$.

b). Suppose $H^s_{ij}(\vec{A})$ has a linear term then
by part a) we have
$s+\llceil\frac{r_{ij}+s}{d-1}\rrceil =1$, which implies $s=0$ and
$r_{ij}\neq 0$. In this case the only solution to (\ref{E:combine})
is $m_{d-r_{ij}}=1$ and $m_\ell=0$ for all
$\ell\neq r_{ij}$.  So the linear monomial is $A_{d-r_{ij}}$.  In the
same vein we obtain the assertion about the constant term.
\end{proof}

For $1\leq n\leq d-1$ and $0\leq t\leq c_n$, let
\begin{eqnarray}
\label{E11}
c_n&:=&\frac{1}{d}
\left(\max_{\sigma\in S_n}\sum_{i=1}^{n}r_{i,\sigma(i)}-
\min_{\sigma\in S_n}\sum_{i=1}^{n}r_{i,\sigma(i)}\right);\\
\label{E10}
S_n^t &:=&\{\sigma\in S_n |
\sum_{i=1}^{n}r_{i,\sigma(i)}
      =\min_{\sigma\in S_n}
       \sum_{i=1}^{n}r_{i,\sigma(i)} +dt\};\\
\label{E12}
f_n^t(\vec{A})&:=&
\sum_{\substack{s_0+s_1+\cdots+s_n=t\\s_0,\ldots,s_n\geq 0}}
\sum_{\sigma\in S_n^{s_0}}\sgn(\sigma)
\prod_{i=1}^{n}H^{s_i}_{i,\sigma(i)}(\vec{A}).
\end{eqnarray}
Note that $c_n\leq n$.
The polynomial $f_n^t(\vec{A})\in\Q[\vec{A}]$
will play a central role in this paper.

\begin{lemma}[Key-Lemma]\label{L:5}
Let $1\leq n\leq d-1$. Then there exists $t$ with
$0\leq t\leq c_n$ such that
the polynomial $f_n^t(\vec{A})\neq 0$.
Let $t_n$ be the least such $t$.
Let $\cV_n$ be the complement in $\A^{d-1}$ of the
variety defined by $f^{t_n}_n=0$.
Then $\cV_n$ is a Zariski dense open subset defined over $\Q$
of $\A^{d-1}$.
\end{lemma}
\begin{proof}
It suffices to prove the first assertion.
  We first show that among the lowest-degree-terms in the formal
  expansion of $\sum_{t=0}^{c_n}f_n^t$ there is a unique
  highest-lexicographic-order-monomial. This suffices because the
  polynomial $f_n^t$ (for some $t$)
whose formal expansion contains this unique
  monomial has to be nonzero.

Partition the summands of the formal expansion
of $\sum_{t=0}^{c_n}f_n^t$ into two parts:
\begin{eqnarray*}
\sum_{t=0}^{c_n}f_n^t
&=&
\sum_{t=0}^{c_n}
{\sum}'\sum_{\sigma\in S_n^{s_0}}\sgn(\sigma)
\prod_{i=1}^{n}H^{s_i}_{i,\sigma(i)}(\vec{A})
+
\sum_{t=0}^{c_n}
{\sum}''\sum_{\sigma\in S_n^{s_0}}\sgn(\sigma)
\prod_{i=1}^{n}H^{s_i}_{i,\sigma(i)}(\vec{A})
\end{eqnarray*}
where $\sum'$
ranges over the set of all $s_0,\ldots,s_n\geq 0$
with $s_1=\cdots=s_n=0$
and $s_0=t$ while $\sum''$ ranges over
the set of all $s_0,\ldots,s_n\geq 0$
with $s_0+\cdots+s_n=t$ and $s_\ell\geq 1$
for some $\ell=1,\ldots,n$.
Denote them by  $H'(\vec{A})$ and $H''(\vec{A})$, respectively.
Note that $S_n=\bigcup_{t=0}^{c_n}S_n^t$, by which we find
\begin{eqnarray*}
H'=\sum_{\sigma\in S_n}\sgn(\sigma)
\prod_{i=1}^{n}H^{0}_{i,\sigma(i)}(\vec{A}).
\end{eqnarray*}
Let $\mu$, $\mu'$ and $\mu''$
denote the lowest degrees in the formal expansions of
$\sum_{t=0}^{c_n}f_n^t$, $H'$ and $H''$, respectively.
By Lemma \ref{L:H}a, we have
$\mu''=
\sum_{i=1}^{n}(s_i+\llceil
\frac{r_{i,\sigma'(i)}+s_i}{d-1}\rrceil)$
for some $\sigma'\in S_n$.
By the definition of $H''$ we have
$s_i\geq 1$ for some $1\leq i \leq n$.
Thus
$\sum_{i=1}^{n}\llceil\frac{r_{i,\sigma'(i)}}{d-1}\rrceil<\mu''$.
On the other hand, we have
\begin{equation*}
\mu\leq\mu'=
\min_{\sigma\in S_n}
      \sum_{i=1}^{n}\llceil\frac{r_{i,\sigma(i)}}{d-1}\rrceil\leq
\sum_{i=1}^{n}\llceil\frac{r_{i,\sigma'(i)}}{d-1}\rrceil.
\end{equation*}
Combining these above, we have $\mu\leq\mu'<\mu''$.
Hence $\mu=\mu'<\mu''$ and it follows that
all degree-$\mu$ monomials
in the formal expansion of
$\sum_{t=0}^{c_n}f_n^t$ lie in the formal expansion of $H'$.

Recall from Lemma \ref{L:H}b that for every
$i$ the lowest-degree-monomial of $H^0_{i,\sigma(i)}$ is $1$ or
$A_{d-r_{i,\sigma(i)}}$ depending on $r_{i,\sigma(i)}=0$ or not,
respectively.
Then the set of degree-$\mu$ monomials of the
formal expansion of $H'$ is equal to the set of
degree-$\mu$ monomials in
$\dot{D}_n|_{A_d=1}$
by a perusal of the definition of $\dot{D}_n$ in (\ref{E:D_n}).
This finishes the proof by Lemma \ref{L:evaluation}.
\end{proof}

\section{Fredholm polygons}
\label{S:4}

Let notations be as in previous sections.
This section will study the shape of Fredholm polygons
of $f\in\A^{d-1}$. We do this by considering
the $p$-adic valuation of the content of
$G_{pi-j}(\vec{A})\in\Z_p[\zeta_p][\vec{A}]$ and that of the
$C_n(\vec{A})\in\Z_p[\zeta_p][[\vec{A}]]$.
We shall consider $G_{pi-j}(\vec{A})$
as formal expressions in $\Z_p[\vec{A}][\gamma]$.

Throughout this section we adopt the following
convention.
Fix an integer $r$
with $1\leq r\leq d-1$ and $\gcd(r,d)=1$.
Let $p$ be a prime that $p\equiv r\bmod d$.
Let $\vec{a}=(a_1,\ldots,a_{d-1})\in (\Q\cap\Z_p)^{d-1}$.
Let $n$ be an integer with $1\leq n\leq d-1$.
For any rational number $R$ let $\gamma^{>R}$ denote the terms in
$\Q_p(\zeta_p)[[\vec{A}]]$ whose coefficients have $p$-adic valuation
$>R/(p-1)$.  We also use it to denote algebraic numbers in
$\Q_p(\zeta_p)$ with $p$-adic valuation $>R/(p-1)$ and this should not
cause any confusion. We define $\gamma^{\geq R}$ analogously.
Let
$$M_n:=\min_{\sigma\in S_n}\sum_{i=1}^{n}\llceil\frac{pi-\sigma(i)}{d}
\rrceil.
$$

\begin{lemma}\label{L:1}
For any $s\geq 0$  we have
\begin{eqnarray}
c_n&=&\max_{\sigma\in S_n}
\sum_{i=1}^{n}\llceil\frac{pi-\sigma(i)}{d}\rrceil-M_n\leq n;
\label{E2}\\
M_n &=&\frac{n(n+1)(p-1)}{2d}+
       \frac{1}{d}\min_{\sigma\in S_n}\sum_{i=1}^{n}r_{i,\sigma(i)};
\label{E1}\\
S^s_n&=&\{\sigma\in S_n|
      \sum_{i=1}^{n}\llceil\frac{pi-\sigma(i)}{d}\rrceil=M_n+s\}.\label{E3}
\end{eqnarray}
\end{lemma}
\begin{proof}
Suppose $\sigma_1,\sigma_2\in S_n$ are minimizer and
maximizer of $\sum_{i=1}^{n}\llceil\frac{pi-\sigma(i)}{d}\rrceil$,
respectively.
Note that $\llceil\frac{pi-j}{d}\rrceil=\frac{pi-j+r_{ij}}{d}$
thus
$$
\max_{\sigma\in S_n}\sum_{i=1}^{n}\llceil\frac{pi-\sigma(i)}{d}\rrceil-M_n
=\frac{1}{d}\left(\max_{\sigma\in S_n}\sum_{i=1}^{n}r_{i,\sigma(i)}
 -\min_{\sigma\in S_n}\sum_{i=1}^{n}r_{i,\sigma(i)}\right)=c_n.
$$
For any $i$, since $1\leq \sigma(i)\leq d-1$, we have
$$
\llceil \frac{pi-\sigma_2(i)}{d} \rrceil
\leq\llceil\frac{pi-\sigma_1(i)}{d}\rrceil+1.
$$
Taking sum both sides and get
\begin{eqnarray*}
\max_{\sigma\in S_n}
\sum_{i=1}^{n}\llceil\frac{pi-\sigma(i)}{d}\rrceil\leq M_n+n.
\end{eqnarray*}
This proves (\ref{E2}).
Since
\begin{eqnarray*}
\sum_{i=1}^{n}\llceil\frac{pi-\sigma(i)}{d}\rrceil
=\frac{n(n+1)(p-1)}{2d}+\frac{1}{d}\sum_{i=1}^{n}r_{i,\sigma(i)},
\end{eqnarray*}
we see that (\ref{E1}) and (\ref{E3}) follows.
\end{proof}

For $0\leq s,t\leq c_n$, and $i,j\geq 1$ let
\begin{eqnarray}\label{E:K}
K^s_{ij}(\vec{A})&:=&
\sum_{\vec{m}\in\cM^s_{ij}}\frac{\vec{A}^{\vec{m}}}
{\vec{m}!(\llceil\frac{pi-j}{d}\rrceil+s-|\vec{m}|)!}.\\
f^t_{n,p}(\vec{A})&:=&
\sum_{\substack{s_0+\cdots+s_n=t\\s_0,\ldots,s_{n}\geq 0}}
\sum_{\sigma\in S_n^{s_0}}\sgn(\sigma)
\prod_{i=1}^{n}K^{s_i}_{i,\sigma(i)}(\vec{A}).
\end{eqnarray}

For $p\geq d^2$, one notes that $H^s_{ij}(\vec{A}),
K^s_{ij}(\vec{A})\in\Z_p[\vec{A}]$, hence
$f_n^t(\vec{A})$, $f^t_{n,p}(\vec{A})\in\Z_p[\vec{A}]$.
But $f_n^t(\vec{A})$ evaluates at $\vec{A}=\vec{a}$ while
$f^t_{n,p}(\vec{A})$ at $\vec{A}=\vec{\hat{a}}$.

\begin{lemma}\label{L:after1}
Let $p\geq (d^2+1)(d-1)$. Then
$f_n^t(\vec{A})\equiv u_nf^t_{n,p}(\vec{A}) \bmod p$
for some $p$-adic unit $u_n$, where the reduction is taken at coefficients.
Moreover, $f_n^t(\vec{a})\equiv u_nf^t_{n,p}(\vec{\hat{a}}) \bmod p$.
\end{lemma}
\begin{proof}
Since $p\geq d^2-1$ we always have
$1\leq \llceil\frac{pi-1}{d}\rrceil+n\leq p-1$.
Then
$$u_n:=\prod_{i=1}^{n}(\llceil\frac{pi-1}{d}\rrceil+n)!$$
is a $p$-adic unit in $\Z_p$.

Recall $\delta{ij}$ defined in (\ref{E:delta}).
It is an elementary exercise to get
\begin{eqnarray*}
\llceil\frac{pi-1}{d}\rrceil&=&\frac{pi+r_{i1}-1}{d}\equiv
\frac{r_{i1}-1}{d}\bmod p\\
\llceil\frac{pi-j}{d}\rrceil &=& \llceil\frac{pi-1}{d}\rrceil-\delta_{ij}
\equiv \frac{r_{i1}-1}{d}-\delta_{ij}\bmod p.
\end{eqnarray*}
Then we have
\begin{eqnarray*}
H^s_{ij}(\vec{A})
&\equiv&
\sum_{\vec{m}\in\cM_{ij}^s}
\frac{(\llceil\frac{pi-1}{d}\rrceil+n)(\llceil\frac{pi-1}{d}\rrceil+n-1)
\cdots(\llceil\frac{pi-j}{d}\rrceil+s+1-|\vec{m}|)}{\vec{m}!}
\vec{A}^{\vec{m}}
\\
&&\equiv
\sum_{\vec{m}\in\cM_{ij}^s}
\frac{(\llceil\frac{pi-1}{d}\rrceil+n) !}
{\vec{m}!(\llceil\frac{pi-j}{d}\rrceil+s-|\vec{m}|)!}
\vec{A}^{\vec{m}}
\\
&&\equiv
(\llceil\frac{pi-1}{d}\rrceil+n)!\ K^s_{ij}(\vec{A})\bmod p.
\end{eqnarray*}
Our first assertion follows easily. The second assertion follows
from the fact that $\vec{a}\equiv \vec{\hat{a}}\bmod p$.
\end{proof}

\begin{proposition}\label{P:11}
Let $p\geq (d^2+1)(d-1)$. For any $1\leq i,j\leq n$
we have
\begin{eqnarray}\label{E:5}
G_{pi-j}(\vec{A})
&=&\sum_{s=0}^{c_n}\gamma^{\llceil\frac{pi-j}{d}\rrceil+s}K^s_{ij}(\vec{A})
         +\gamma^{>\llceil\frac{pi-j}{d}\rrceil+c_n}.\\
\label{E:6}
\det\{G_{pi-j}(\vec{A})\}_{1\leq i,j\leq n}
&=&\sum_{t=0}^{c_n}
\gamma^{M_n+t}f^t_{n,p}(\vec{A})
+\gamma^{>M_n+c_n}.
\end{eqnarray}
\end{proposition}
\begin{proof}
For $0\leq s\leq c_n$, $m_\ell\geq 0$
and $m_1+\cdots+m_d=\llceil\frac{pi-j}{d}\rrceil+s$,
since $p\geq d^2-1$, we have
$m_\ell\leq \llceil\frac{pi-j}{d}\rrceil+c_n\leq p-1$.
And by (\ref{E:33}) and (\ref{E:3}), we have
\begin{eqnarray*}
G_{pi-j}(\vec{A})
&=&\sum_{
\substack{m_1+\cdots+m_d\leq \llceil\frac{pi-j}{d}\rrceil+c_n\\
\sum_{\ell=1}^{d}\ell m_\ell=pi-j}}
\lambda_{m_1}\cdots \lambda_{m_d}\vec{A}^{\vec{m}}
+\gamma^{>\llceil\frac{pi-j}{d}\rrceil+c_n}\\
&=&\sum_{s=0}^{c_n}\sum
\frac{\gamma^{m_1+\cdots+m_d}\vec{A}^{\vec{m}}}{m_1!\cdots m_d!}
 +\gamma^{>\llceil\frac{pi-j}{d}\rrceil+c_n}
\end{eqnarray*}
where the last sum ranges over all $m_\ell\geq 0$ such that $m_1+\cdots
+m_d=\llceil\frac{pi-j}{d}\rrceil+s$ and $\sum_{\ell=1}^{d}\ell
m_\ell=pi-j$.  It is easy to see that this is a subset of
$\cM^s_{ij}$.  Conversely, if $\vec{m}\in \cM^s_{ij}$ then
$$
\sum_{\ell=1}^{d-1}m_\ell\leq \sum_{\ell=1}^{d-1}\ell
m_{d-\ell}=r_{ij}+ds\leq d-1+ds
\leq
\llceil\frac{pi-j}{d}\rrceil+s
$$
since $p\geq (d^2+1)(d-1)$.
Set $m_d=\llceil\frac{pi-j}{d}\rrceil+s-\sum_{\ell=1}^{d-1}m_\ell$,
then $m_1+\cdots+m_{d}=\llceil\frac{pi-j}{d}\rrceil+s$ and
$\sum_{\ell=1}^{d}\ell m_\ell=pi-j$ where $m_\ell\geq 0$.
Thus we have
\begin{eqnarray*}
G_{pi-j}(\vec{A})
&=&\sum_{s=0}^{c_n}
   \gamma^{\llceil\frac{pi-j}{d}\rrceil+s}
   \sum_{\vec{m}\in\cM^s_{ij}}\frac{\vec{A}^{\vec{m}}}
   {\vec{m}!(\llceil\frac{pi-j}{d}\rrceil+s-|\vec{m}|)!}
   +\gamma^{>\llceil\frac{pi-j}{d}\rrceil+c_n}.
\end{eqnarray*}

To prove (\ref{E:6}) we have
\begin{eqnarray*}
\det\{G_{pi-j}(\vec{A})\}_{1\leq i,j\leq n}
&=&\sum_{\sigma\in S_n}\sgn(\sigma)
\prod_{i=1}^{n}G_{pi-\sigma(i)}(\vec{A})\\
&=&\sum_{\sigma\in S_n}\sgn(\sigma)
\prod_{i=1}^{n}\sum_{s_i=0}^{c_n}\left(
\gamma^{\llceil\frac{pi-\sigma(i)}{d}\rrceil+s_i}K^{s_i}_{i,\sigma(i)}(\vec{A})
+\gamma^{>\llceil\frac{pi-\sigma(i)}{d}\rrceil+c_n}\right)
\\
&=&\sum_{s_0=0}^{c_n}\sum_{\sigma\in S_n^{s_0}}\sgn(\sigma)
   \sum_{\ell=0}^{c_n-s_0}\gamma^{M_n+s_0+\ell}\sum_{s_1+\cdots+s_n=\ell}
   \prod_{i=1}^{n}K^{s_i}_{i,\sigma(i)}(\vec{A})+ \gamma^{>M_n+c_n}\\
&=&\sum_{t=0}^{c_n}\gamma^{M_n+t}
\left(\sum_{s_0+\cdots+s_n=t}\sum_{\sigma\in S_n^{s_0}}\sgn(\sigma)
\prod_{i=1}^{n}K^{s_i}_{i,\sigma(i)}(\vec{A})\right)+\gamma^{>M_n+c_n},
\end{eqnarray*}
where the second equality follows from (\ref{E:5})
and the third from Lemma \ref{L:1}.
\end{proof}

\begin{lemma}\label{L:2}
Let $p\geq (d^2+1)(d-1)$.  Then
  $\ord_pC_n(\vec{\hat{a}})\geq\frac{M_n+t_n}{p-1}$ for all
$\vec{a}\in (\Z_p\cap\Q)^{d-1}$,
and the equality
  holds if and only if $\vec{\bar{a}}\in\cV_n(\F_p)$.
\end{lemma}
\begin{proof}
First we show that
\begin{eqnarray}\label{E:C_n2}
C_n(\vec{A})&=&\sum_{t=0}^{c_n}
\gamma^{M_n+t}f^t_{n,p}(\vec{A}) +\gamma^{>M_n+c_n}.
\end{eqnarray}
By (\ref{E:1}) and (\ref{E:6})
it suffices to show that if there is a $t$ with $u_t>n$
then
\begin{eqnarray}\label{E:contradiction}
\min_{\sigma\in S_n}\ord_p \prod_{t=1}^{n} G_{pu_t-u_{\sigma(t)}}(\vec{A})
&>& \frac{M_n+c_n}{p-1}.
\end{eqnarray}
Since $\sum_{t=1}^{n}u_t>\sum_{i=1}^{n}i$, we have
\begin{eqnarray}\label{E:twin1}
\frac{1}{p-1}\sum_{t=1}^{n}\llceil\frac{pu_t-u_{\sigma(t)}}{d}\rrceil
&\geq & \frac{1}{p-1}\sum_{t=1}^{n}\frac{pu_t-u_{\sigma(t)}}{d}\nonumber\\
&=&\frac{1}{d}\sum_{t=1}^{n}u_t\nonumber\\
&\geq&\frac{1}{d}\sum_{i=1}^{n}i+\frac{1}{d}
   \nonumber\\
&=& \frac{1}{p-1}\sum_{i=1}^{n}\frac{pi-\delta(i)}{d}+\frac{1}{d}
\end{eqnarray}
for any $\delta\in S_n$.
For $p\geq (d^2+1)(d-1)> d^2-d+1$ we have
\begin{eqnarray}\label{E:twin2}
\frac{1}{p-1}\sum_{i=1}^{n}\llceil\frac{pi-\delta(i)}{d}\rrceil
&\leq&\frac{1}{p-1}\sum_{i=1}^{n}\frac{pi-\delta(i)}{d}+\frac{n}{p-1}
  \nonumber\\
&\leq&\frac{1}{p-1}\sum_{i=1}^{n}\frac{pi-\delta(i)}{d}+\frac{d-1}{p-1}
  \nonumber\\
&<&\frac{1}{p-1}\sum_{i=1}^{n}\frac{pi-\delta(i)}{d}+\frac{1}{d}
\end{eqnarray}
for any $\delta\in S_n$.
Therefore,
\begin{eqnarray*}
\min_{\sigma\in S_n}\ord_p \prod_{t=1}^{n}G_{pu_t-u_{\sigma(t)}}(\vec{A})
&\geq&\frac{1}{p-1}\min_{\sigma\in S_n}
\sum_{t=1}^{n}\llceil\frac{pu_t-u_{\sigma(t)}}{d}\rrceil\\
&>&\frac{1}{p-1}\max_{\delta\in S_n}
\sum_{i=1}^{n}\llceil\frac{pi-\delta(i)}{d}\rrceil
=\frac{M_n+c_n}{p-1}
\end{eqnarray*}
where the first inequality is due to (\ref{E:Bound}), the second
inequality
by (\ref{E:twin1}) and (\ref{E:twin2}), and the last by (\ref{E2}).

Let $0\leq t< t_n$. We have $f^t_n(\vec{A})=0$ and hence
by Lemma \ref{L:after1} we have $f^t_{n,p}(\vec{\hat{a}})
\equiv 0\bmod p$.
So
$$\ord_p(\gamma^{M_n+t}f^t_{n,p}(\vec{\hat{a}}))
\geq \frac{M_n+t}{p-1}+1>\frac{M_n+c_n}{p-1}.
$$
Therefore,
for all $\vec{a}\in(\Z_p\cap\Q)^{d-1}$
by (\ref{E:C_n2} we have
$$
C_n(\vec{\hat{a}})=f_{n,p}^{t_n}(\vec{\hat{a}})\gamma^{M_n+t_n}
+\gamma^{>M_n+t_n}.
$$
So
$$
\ord_pC_n(\vec{\hat{a}})\geq \frac{M_n+t_n}{p-1}
$$
and the equality holds if and only if
$f_{n,p}^{t_n}(\vec{\hat{a}})\equiv f_n^{t_n}(\vec{a})
\not\equiv 0\bmod p$.
This proves the lemma.
\end{proof}

\section{Generic Newton polygons}\label{S:Proof}

Let the {\em generic Newton polygon} of $\A^d$ over $\F_p$
be the lowest Newton polygon over all $\bar{f}\in\A^d(\F_p)$,
that is,
$$\GNP(\A^d;\F_p):=\inf_{\bar{f}\in\A^d(\F_p)}\NP(\bar{f}).$$
Note that it is equal to
$\inf_{f\in\A^d(\Z_p\cap\Q)}\NP(f\otimes\F_p)$
and one does not know {\it a priori} whether this infimum exists.
Note that Wan has shown that the generic Newton polygon over
$\bar\F_p$ defined by
$\GNP(\A^d;\bar\F_p)
:=\inf_{\bar{f}\in\A^d(\bar\F_p)}\NP(\bar{f})$ exists
by Grothendieck specialization theorem
(see \cite[Section 1.1]{Wan:2}). In the theorem below
we show that $\GNP(\A^d;\F_p)$ exists
for $p$ large enough. One may ask if it is true that
$\GNP(\A^d;\F_p)=\GNP(\A^d;\bar\F_p)$
for $p$ large enough.

We shall proceed to prove Theorem \ref{T:2} below by first
introducing some notations.
Let $\epsilon_0=0$ and for $1\leq n\leq d-1$ let
\begin{eqnarray}\label{E:epsilon}
\epsilon_n&:=&
\frac{\min_{\sigma\in S_n}\sum_{i=1}^{n}r_{i,\sigma(i)}+dt_n}{d(p-1)},
\end{eqnarray}
where $r_{ij}$ and $t_n$ are defined in section \ref{S:31}
and Lemma \ref{L:5}, respectively.
One observes easily
\begin{eqnarray}\label{E:M}
\frac{M_n+t_n}{p-1}&=&
\frac{n(n+1)}{2d}+\epsilon_n.
\end{eqnarray}
Note that $0\leq r_{ij}\leq d-1$ for all $1\leq i,j\leq d-1$,
and $t_n\leq c_n\leq n\leq d-1$ by (\ref{E2}),
so we have $\epsilon_n\leq \frac{n(2d-1)}{d(p-1)}$. Thus
$\epsilon_n$ goes to $0$ as $p$ approaches $\infty$.

For every integer $r$ with $1\leq r\leq d-1$ and $\gcd(r,d)=1$, let
$\cW_r:=\bigcap_{n=1}^{d-1}\cV_n$ (recall from
Key-Lemma \ref{L:5} that $\cV_n$ consists of all $f\in\A^{d-1}$
whose coefficients $\vec{a}$ satisfy
$f_n^{t_n}(\vec{a})\neq 0$).
Let $\cW:=\bigcap_{\substack{1\leq
    r\leq d-1\\\gcd(r,d)=1}}\cW_r$.  Consider the natural projection
map $\iota:\A^d\rightarrow \A^{d-1}$ by
$\iota(f)=\vec{a}=(a_1,\ldots,a_{d-1})$ for
every $f=x^d+a_{d-1}x^{d-1}+\cdots+a_0\in\A^d$. Let
$\cU:=\iota^{-1}(\cW)$.
For every residue class $r$
denote by $f_{n,r}^{t_n}$ the
$f_n^{t_n}$ in Lemma \ref{L:5}, then
$\cU$ consists of all $f\in\A^d$ whose coefficients satisfy
$\prod_{r}\prod_{n=1}^{d-1}f_{n,r}^{t_n}(\vec{a})\neq 0$
where $r$ ranges over all $1\leq r\leq d-1$ coprime to $d$.
Since $\prod_{r}\prod_{n=1}^{d-1}f_{n,r}^{t_n}$
is a nonzero polynomial over $\Q$ by
Lemma \ref{L:5}, one concludes that
$\cU$ is Zariski dense open in $\A^d$ over $\Q$.
One notes that, even though $\cU(\F_p)$ is not necessarily nonempty,
it is nonempty when $p$ is large enough.

\begin{theorem}\label{T:2}
Let notations be as above.\\
a) For $p$ large enough (depending only on $d$) $\GNP(\A^d;\F_p)$ exists and
is equal to the lower convex hull of
points $(n, \frac{n(n+1)}{2d}+\epsilon_n)$ for $0\leq n\leq d-1$,
each of which is a vertex.\\
b) Fix $f\in\A^d(\Q)$.
For $p$ large enough (depending only on $d$ and $f$) we have
\begin{eqnarray*}
\NP(f\otimes\F_p)\geq\GNP(\A^d;\F_p)
\end{eqnarray*}
where the equality holds for all $p$ large enough if and only if
$f\in\cU(\Q)$. Here $\geq$ means ``lies above''.\\
c) For $f\in\cU(\Q)$ we have
$$
\lim_{p\rightarrow\infty}\NP(f\otimes\F_p) = \HP(\A^d).
$$
\end{theorem}
\begin{proof}
a) Because of (\ref{E:98}), we consider $f(x)\in\A^d(\Z_p\cap\Q)$
with no constant term, that is, $f(x)=x^d+\sum_{i=1}^{d-1}a_ix^i$.
By Lemma \ref{L:2}, for $p$ large enough we have
$$\ord_pC_n(\vec{\hat{a}})\geq \frac{n(n+1)}{2d}+\epsilon_n$$
for all
$0\leq n\leq d-1$ and the equality holds if and only if
$\vec{\bar{a}}\in\cW(\F_p)$.
On the other hand, by the remarks preceding the theorem,
$\epsilon_n$ approaches $0$.
Thus for $p$ large enough
the lower convex hull of points
$(n,\frac{n(n+1)}{2d}+\epsilon_n)$ with $0\leq n\leq d-1$
passes all these points as vertices.
By Proposition \ref{P:Dwork}, for $p$ large enough,
\begin{eqnarray}\label{E:bn}
\ord_pb_n(\vec{a})\geq \frac{n(n+1)}{2d}+\epsilon_n
\end{eqnarray}
and the equality holds if and only if
$\vec{\bar{a}}\in\cW(\F_p)$.
Now a) clearly follows.

b) Now let $f(x)=x^d+a_{d-1}x^{d-1}+\cdots+a_1x+a_0\in\A^d(\Q)$.
Note that (\ref{E:bn}) says for $p$ large enough,
\begin{eqnarray*}
\NP(f\otimes\F_p)=\NP((x^d+a_{d-1}x^{d-1}+\cdots+a_1x)\otimes\F_p)
\geq\GNP(\A^d;\F_p)
\end{eqnarray*}
where the equality holds if and only if
$(\bar{a}_1,\ldots,\bar{a}_{d-1})\in\cW(\F_p)$.
Note that a rational number $N$ is nonzero if and only if
$N$ is not divisible by all primes large enough.
Thus for $p$ large enough the above equality holds if and only if
$(a_1,\cdots,a_{d-1})\in\cW(\Q)$, that is, $f\in\cU(\Q)$.
This proves b). Note that c) follows from a) and b).
\end{proof}

\begin{remark}
1)
Let $d\geq 3$.
Let {\em generic polynomial}
$F_d:=\prod_{r}\prod_{n=1}^{d-1}f_{n,r}^{t_n}$
where $r$ ranges over $1\leq r\leq d-1$ coprime to $d$.
From the theorem above, the set of
polynomials $f(x)=x^{d}+\cdots+a_1x+a_0\in\A^d(\Q)$ with
$\NP(f\otimes\F_p)=\GNP(\A^d)$ corresponds
precisely to the set of $(a_0,\ldots,a_{d-1})\in\Q^d$ with
$F_d|_{\vec{A}=(a_1,\ldots,a_{d-1})}\neq 0$.

2) In practice, for any $d\geq 3$
one may compute the polynomial
$P_d:=\prod_{r}\prod_{n=1}^{\lceil{\frac{d-1}{2}}\rceil}f_n^{t_n}$ in
$\Q[\vec{A}]$ where $r$ ranges over all $2\leq r\leq d-1$ with
$\gcd(r,d)=1$. (Remark:
One notes that the $r=1$ case is explained
in remarks before Theorem \ref{T:1}. One also notes that
$\NP(f\otimes\F_p)$ is symmetric in the sense that
every slope $\alpha$ segment comes with a slope $1-\alpha$ segment
with the same length.)
Then every $f(x)=x^d+a_{d-1}x^{d-1}+\cdots+a_0\in\A^d(\Q)$
with $P_d|_{\vec{A}=(a_1,\ldots,a_{d-1})}\neq 0$
satisfies $\lim_{p\rightarrow\infty}\NP(f\otimes\F_p)=\HP(\A^d).$
\end{remark}

\section{Generic Newton polygon for $x^d+ax$}
\label{S:6}

Recall $d\geq 3$.
In this section we consider the Newton polygon of
the $L$ function of exponential sums of
$f(x)=x^d+ax$ over $\Q$.
This family has drawn some attentions
recently (see \cite{Yang} for some progress).
When
$a=0$ see Remark \ref{R:13}b.
Let $\A^d(1)$ denote the space of all such $f(x)$ with parameter $a$.
Let $\GNP(\A^d(1);\F_p)$ be the corresponding analog of $\GNP(\A^d;\F_p)$.

Let $r$ be $1\leq r\leq d-1$ coprime to $d$.
Recall that $r'_{ij}$ is the least nonnegative residue of $ri-j$ mod $d$.
That is, $r'_{ij}=ri-j-d\llfloor\frac{ri-j}{d}\rrfloor$.
Let $M'_n:=\min_{\sigma\in S_n}\sum_{i=1}^{n}r'_{i,\sigma(i)}$.
Let $S'_n$ be the subset of $\sigma\in S_n$
with $\sum_{i=1}^{n}r'_{i,\sigma(i)}=M'_n$.
Let $\epsilon'_0:=0$; for $n\geq 1$
and for $p\equiv r\bmod d$ let
$$\epsilon'_n:=\frac{(d-1)M'_n}{d(p-1)}.$$

\begin{lemma}\label{L:Sprime}
Let $1\leq n\leq d-1$ and $p \equiv r \bmod d$.
The following statements are equivalent:
\begin{enumerate}
\item[1)] $\sigma\in S'_n$;
\item[2)] $\sigma(i)\leq r'_{i1}+1$ for all $1\leq i\leq n$;
\item[3)] $r'_{i,\sigma(i)}=r'_{i1}-\sigma(i)+1$ for all $1\leq i\leq n$;
\item[4)] $\llfloor\frac{pi-1}{d}\rrfloor=\llfloor\frac{pi-\sigma(i)}{d}\rrfloor$
for all $1\leq i\leq n$.
\end{enumerate}
\end{lemma}
\begin{proof}
Define $\delta'_{ij}:=0$ if $j\leq r'_{i1}+1$ and
$\delta'_{ij}:=1$ if $j>r'_{i1}+1$.
From Lemma \ref{L:r_n} one notes that
$r'_{11}+1, \ldots, r'_{n1}+1$ are
$n$ distinct integers in the interval $[1,d-1]$.
So there exists $\sigma\in S_n$ such that
$\sigma(i)\leq r'_{i1}+1$ for every $1\leq i\leq n$,
that is, $\delta'_{i,\sigma(i)}=0$ for every $1\leq i\leq n$.
Thus $\min_{\sigma\in S_n}\sum_{i=1}^{n}\delta'_{i,\sigma(i)}=0$
and it is achieved if and only if 2) holds.

By recalling Lemma \ref{L:r_n},
it is straightforward to see that
\begin{eqnarray*}\label{E:rprime}
r'_{ij} &=& r'_{i1} - j  + 1 +\delta'_{ij}(d-1).
\end{eqnarray*}
Thus for any $\sigma\in S_n$,
$$\sum_{i=1}^{n}r'_{i,\sigma(i)}=
\sum_{i=1}^{n}(r'_{i1}-\sigma(i)+1)+(d-1)\sum_{i=1}^{n}\delta'_{i,\sigma(i)}
=\sum_{i=1}^{n}r'_{i1} -\frac{n(n-1)}{2}
+(d-1)\sum_{i=1}^{n}\delta'_{i,\sigma(i)}
$$
One notes that  1) holds if and only if
$\sum_{i=1}^{n}r'_{i,\sigma(i)}$ achieves its minimum,
and if and only if $\sum_{i=1}^{n}\delta'_{i,\sigma(i)}=0$
by the previous paragraph.
Thus 1), 2) and 3) are equivalent to each other.
Since
$r'_{ij}=pi-j-d\lfloor\frac{pi-j}{d}\rfloor$,
it is easy to see 3) and 4) are equivalent to each other.
This proves the lemma.
\end{proof}

By the lemma above,
$M'_n=\sum_{i=1}^{n}(r'_{i1}-\sigma(i)+1)$.
So one gets an explicit formula
$$\epsilon'_n =
\frac{(d-1)(\sum_{i=1}^{n}r'_{i1} - \frac{n(n-1)}{2})}{d(p-1)}.
$$
Note that $\epsilon'_n$ converges to $0$ as $p$ approaches $\infty$.

\begin{theorem}\label{T:3}
a) For $p$ large enough (depending only on $d$) $\GNP(\A^d(1);\F_p)$
exists and is equal to the lower convex hull of
points $(n, \frac{n(n+1)}{2d}+\epsilon'_n)$ for $0\leq n\leq d-1$,
each of which is a vertex.\\
b) Fix $f=x^d+ax\in\A^d(\Q)$.
For $p$ large enough (depending only on $d$ and $a$) we have
\begin{eqnarray*}
\NP(f\otimes\F_p)\geq\GNP(\A^d(1);\F_p),
\end{eqnarray*}
where the equality holds for all $p$ large enough if and only if
$a\neq 0$.
Here $\geq$ means ``lies above''.\\
c) For any $a\neq 0$ we have
$$
\lim_{p\rightarrow\infty}\NP((x^d+ax)\otimes\F_p) = \HP(\A^d).
$$
\end{theorem}

\begin{lemma}\label{L:A1}
Let $p\equiv r\bmod d$.
Let $a\in\Q\cap\Z_p$ and let
$\hat{a}$ be the Teichm\"uller lifting of $a\bmod p$.
Let $p\geq d$. For any $1\leq i,j\leq d-1$ we have
$$G_{pi-j}= \gamma^{r'_{ij}+\llfloor\frac{pi-j}{d}\rrfloor}
\hat{a}^{r'_{ij}}\frac{1}{r'_{ij}!\llfloor\frac{pi-j}{d}\rrfloor
  !}+\gamma^{>r'_{ij}+\llfloor\frac{pi-j}{d}\rrfloor}.$$
\end{lemma}
\begin{proof}
Note that
$
G_{pi-j}
=\sum \lambda_{m_1}\lambda_{m_d}\hat{a}^{m_1}
$
where the sum ranges in $m_1+dm_d=pi-j$ with $m_1,m_d\geq 0$.
But in this range of $m_1$ and $m_d$, one notices that
the minimum of
$m_1+m_d$ is achieved precisely at $m_1=r'_{ij}$
and $m_d=\lfloor\frac{pi-j}{d}\rfloor$, that is,
$\min(m_1+m_d)=r'_{ij}+\lfloor\frac{pi-j}{d}\rfloor$.
The rest of the proof is analogous to Proposition \ref{P:11}.
\end{proof}

\begin{lemma}\label{L:A2}
Let $p\equiv r\bmod d$ and $p\geq (d-1)^3+2$.
Let $1\leq n\leq d-1$.
Then we have
$$C_n=\gamma^{(p-1)(\frac{n(n+1)}{2d}+\epsilon'_n)}
\hat{a}^{M'_n}f'_{n,p}+\gamma^{>(p-1)(\frac{n(n+1)}{2d}+\epsilon'_n)},$$
where
$$
f'_{n,p}:=\sum_{\sigma\in S'_n}
\sgn(\sigma)
\prod_{i=1}^{n}\frac{1}{r'_{i,\sigma(i)}!
\llfloor\frac{pi-\sigma(i)}{d}\rrfloor !}.
$$
\end{lemma}
\begin{proof}
The proof is analogous to Lemma \ref{L:2}, so we will
only give an outline.
First one shows that for $1\leq n\leq d-1$ one has
\begin{eqnarray}\label{E:A2}
C_n&=&\sum_{\sigma\in S_n}\sgn(\sigma)\prod_{i=1}^{n}G_{pi-\sigma(i)}
+\gamma^{\geq (p-1)(\frac{n(n+1)}{2d}+\frac{1}{d})}.
\end{eqnarray}
Since $r'_{ij}=pi-j-d\llfloor\frac{pi-j}{d}\rrfloor$,
we have
$$
r'_{ij}+\llfloor\frac{pi-j}{d}\rrfloor =
\frac{pi-j}{d}+\frac{d-1}{d}r'_{ij}.
$$
Thus
\begin{eqnarray*}
\min_{\sigma\in S_n}\sum_{i=1}^{n}
\left(r'_{i,\sigma(i)}+\llfloor\frac{pi-\sigma(i)}{d}\rrfloor\right)
&=&\frac{(p-1)n(n+1)}{2d}+\frac{(d-1)M'_n}{d}\\
&=& (p-1)(\frac{n(n+1)}{2d}+\epsilon'_n).
\end{eqnarray*}
Consequently the minimum is achieved precisely at all
$\sigma\in S'_n$.
Note that
$p\geq (d-1)^3+2$ implies that
$(p-1)(\frac{n(n+1)}{2d}+\frac{1}{d})>
(p-1)(\frac{n(n+1)}{2d}+\epsilon'_n)$.
By (\ref{E:A2}) and Lemma \ref{L:A1} we have
\begin{eqnarray*}
C_n&=&\gamma^{(p-1)(\frac{n(n+1)}{2d}+\epsilon'_n)}
\hat{a}^{M'_n}
\sum_{\sigma\in S'_n}
\sgn(\sigma)
\prod_{i=1}^{n}\frac{1}{r'_{i,\sigma(i)}!
\llfloor\frac{pi-\sigma(i)}{d}\rrfloor !}\\
&&+\gamma^{>(p-1)(\frac{n(n+1)}{2d}+\epsilon'_n)}.
\end{eqnarray*}
The lemma follows.
\end{proof}

\begin{lemma}\label{L:A3}
Let notation and hypothesis be as in Lemma \ref{L:A2}.
Then $$\ord_pC_n\geq \frac{n(n+1)}{2d}+\epsilon'_n$$
and the equality holds if and only if $a\not\equiv 0\bmod p$.
\end{lemma}
\begin{proof}
Let
\begin{eqnarray*}
u_n&:=&\prod_{i=1}^{n}r'_{i1}!(\llfloor\frac{pi-1}{d}\rrfloor !).
\end{eqnarray*}
By Lemma \ref{L:Sprime}, one sees that
\begin{eqnarray*}
u_nf'_{n,p}
&=&
\sum_{\sigma\in S'_n}\sgn(\sigma)\prod_{i=1}^{n}\frac{r'_{i1}!}
{r'_{i,\sigma(i)}!}.
\end{eqnarray*}
By Lemma \ref{L:Sprime},
we have
\begin{eqnarray*}
u_nf'_{n,p}
&=&
\sum_{\sigma\in S'_n}\sgn(\sigma)\prod_{i=1}^{n}
(r'_{i1}(r'_{i1}-1)\cdots (r'_{i1}-(\sigma(i)-2))\\
&=&
\sum_{\sigma\in S_n}\sgn(\sigma)\prod_{i=1}^{n}
(r'_{i1}(r'_{i1}-1)\cdots (r'_{i1}-(\sigma(i)-2))
\end{eqnarray*}
where we set
$(r'_{i1}(r'_{i1}-1)\cdots (r'_{i1}-(\sigma(i)-2)):=1$ if
$\sigma(i)=1$.
One observes that this is equal to the determinant of a matrix
$M$ shown as below
\begin{eqnarray*}
M
&=&
\left[
\begin{array}{cccc}
1&r'_{11}&r'_{11}(r'_{11}-1)&\cdots\\
1&r'_{21}&r'_{21}(r'_{21}-1)&\cdots\\
 &      &    \vdots      &      \\
1&r'_{n1}&r'_{n1}(r'_{n1}-1)&\cdots
\end{array}
\right].
\end{eqnarray*}
Under natural column transformation $M$
becomes a Vandermonde matrix, we get
\begin{eqnarray*}
u_nf'_{n,p}=\det M
&=& \det\left[
\begin{array}{cccc}
1&r'_{11}&(r'_{11})^2&\cdots\\
1&r'_{21}&(r'_{21})^2&\cdots\\
 &      &    \vdots      &      \\
1&r'_{n1}&(r'_{n1})^2&\cdots
\end{array}
\right]= \prod_{1\leq i<k\leq n}(r'_{k1}-r'_{i1}).
\end{eqnarray*}
Similar as in Lemma \ref{L:r_n}, one notes that
$r'_{i1}\neq r'_{k1}$ for any $i<k$.
One also notes that $u_n$ is a $p$-adic unit.
Therefore, $f'_{n,p}\not\equiv 0\bmod p$ for all $p$.
\end{proof}

\begin{proof}[Proof of Theorem \ref{T:3}]
Theorem \ref{T:3} follows from Lemma \ref{L:A3},
using the same arguments as in the proof of
Theorem \ref{T:2}.
\end{proof}

\end{document}